\documentclass{elsarticle}

\usepackage[utf8]{inputenc}
\usepackage[]{todonotes}
\usepackage[ruled]{algorithm2e}
\usepackage{endnotes}
\usepackage{pgfplots}
\usepackage{lscape}
\usepackage{setspace}
\usepackage{url}
\usepackage{adjustbox}
\pgfplotsset{compat=newest}

\usepackage{verbatim}

\usepackage{amsmath,amsfonts,amsthm,amssymb, mathrsfs}
\def\setR{\mathbb{R}}
\def\setZ{\mathbb{Z}}

\def\ds{\displaystyle}

\usepackage{footnote}
\makesavenoteenv{tabular}

\graphicspath{{./figures/}}

\usepackage{comment}
\usepackage{graphics}
\usepackage{graphicx}
\usepackage{adjustbox}
\usepackage{tikz}
\usepgfplotslibrary{dateplot}
\usepgfplotslibrary{groupplots}
\usepackage{caption}
\usepackage{subcaption}
\usepackage{longtable}
\usepackage{psfrag}
\usepackage{soul}
\usepackage{mwe}

\usepgfplotslibrary{fillbetween}

\usetikzlibrary{er,positioning, calc, shapes, arrows, chains, fit}

\definecolor{myblue}{RGB}{19,145,215}
\definecolor{mygreen}{RGB}{80,176,50}
\definecolor{myred}{RGB}{222,36,16}

\usepackage{lipsum}
\makeatletter
\def\ps@pprintTitle{%
 \let\@oddhead\@empty
 \let\@evenhead\@empty
 \def\@oddfoot{}%
 \let\@evenfoot\@oddfoot}
\makeatother

\begin{document}
  \begin{frontmatter}

\title{\textbf{PUSH}: a primal heuristic based on Feasibility \textbf{PU}mp and \textbf{SH}ifting} 

\author[1]{Valerio Agasucci\footnote{PhD student. }$^{, }$}
\ead{agasucci@diag.uniroma1.it}

\author[1]{Corrado Coppola$^{1, }$}
\ead{corrado.coppola@uniroma1.it}

\author[2]{Giorgio Grani\footnote{Research fellow, Corresponding author. }$^{, }$}
\ead{g.grani@uniroma1.it}

\address[1]{Sapienza University of Rome, Dep. of Computer Science, Control and Management Engineering, Rome, Italy
}

\address[2]{Sapienza University of Rome, Dep. of Statistical Sciences, Rome, Italy
}

\begin{abstract}
This work describes PUSH, a primal heuristic combining Feasibility Pump and Shifting. The main idea is to replace the rounding phase of the Feasibility Pump with a suitable adaptation of the Shifting and other rounding heuristics. The algorithm presents different strategies, depending on the nature of the partial rounding obtained. In particular, we distinguish when the partial solution is feasible, infeasible with potential candidates and infeasible without candidates. We used a threshold to indicate the percentage of variables to round with our algorithm and which other to round to the nearest integer. Most importantly, our algorithm tackles directly equality constraints without duplicating rows.

We select the parameters of our algorithm on the 19 instances provided for the Mip Competition 2022.
Finally, we compared our approach to other start heuristics, like Simple Rounding, Rounding, Shifting and Feasibility Pump on the first 800 MIPLIB2017 instances ordered by the number of non zeros. 

\end{abstract}
\end{frontmatter}
\parindent 0cm

\section{ Introduction}
\label{sec:Introduction}
Mixed Integer Programming (MIP) is a powerful tool to model a vast amount of scientific, industrial, and economical problems, covering the most disparate fields, such as logistics, assignment, and capital management. Since MIP linear problems include both continuous and integer variables, they are proved to belong to the NP-hard class (see \cite{schrijver2003combinatorial} for a more detailed analysis), meaning that they are not solvable in polynomial time. The complete exploration of the integer feasible set, whose cardinality grows exponentially with the number of variables, is yet possible to achieve the optimal solution, but for most of the practically significant instances, it would require unacceptable computational effort. In fact, the only way to solve to optimality any mixed-integer problem is to apply some of the well-known Branch and Bound techniques. However, despite combinatorial optimization community provided a great deal of these algorithms, for which the reader should refer to \cite{lawler1966branch,mitten1970branch,boyd2007branch}, MIP problems complexity is inherent with their belonging to NP-hard class. Therefore, when tackling MIP problems, one either seeks particular structures allowing to bring down the complexity, such as the availability, for a given class of problems, of the optimal formulation or exploits cutting plane generation to dramatically reduce the feasible region dimension. However, we often encounter MIP problems without having any prior knowledge of possible structures and, thus, pursuing the globally optimal solution could be in practice impossible or inefficient, since for our purpose a sub-optimal approximation is considered to be good enough. This makes heuristics one of the most widespread and feasible ways to achieve sub-optimal solutions of MIP problems within an affordable computational time.

For the purpose of highlighting the perspective of our research, we can define two classes of MIP heuristics: improvement heuristics and start heuristics. 

Improvement heuristics try to improve a solution from a given feasible one. One of the most widespread and successful is Local Branching, developed by Fischetti et al. in \cite{fischetti2003local}, which is based on the feasible set reduction employing local branching cuts, i.e., adding specific constraints to eliminate those feasible regions, where we do not expect to find promising solutions. Local Branching has been recently enhanced in \cite{hansen2006variable} and extended in \cite{fischetti2008repairing}, where it is used to retrieve MIP-feasibility. Relaxation Induced Neighborhood Search (RINS) \cite{danna2005exploring} is another example of improvement heuristic, based on the exploration of some promising neighborhoods that are constructed with local search techniques. RINS has also been successfully implemented in the framework of a genetic algorithm in  \cite{momeni2016genetic}. A similar approach is used by Ghosh in \cite{ghosh2007dins}, where Distance Induced Neighborhood Search (DINS) is proposed. However, alongside more traditional heuristics, novel paradigms are gaining increasing interest in the research community, such as evolutionary algorithms, proposed by Rothberg in  \cite{rothberg2007evolutionary} and successfully applied to solve industrial problems in \cite{hung1999evolutionary,medaglia2009hybrid}, as well as more recent artificial intelligence-based methods to apply the Branch and Bound paradigm \cite{nair2020solving}.

Start heuristics, which are the main object of our research, try to find a first integer feasible solution, given any LP-feasible point, which is usually the optimal solution of the LP-relaxation. The main task of a start heuristic is to find a first feasible point in the fastest time possible, which often results in a challenging problem, especially when only using the LP-toolkit to speed up the algorithm. Several different approaches have been proposed in the literature, and, despite a thorough analysis of the state-of-the-art is beyond the limits of this research, we briefly introduce the most widespread and promising ones.

The first class of start heuristics exploits the most trivial operation that can be performed when dealing with an initial LP-feasible but not integer point, the rounding of fractional variables. The most basic example of this type of heuristics is the Simple Rounding, introduced by Balas in \cite{balas1980pivot} only for binary problems; in this case, only variables that can be rounded without losing the LP-feasibility are updated, and, whenever such variables do not exist, the algorithm returns a failure outcome. In \cite{balas1986pivot} the algorithm was extended to general MIP problems and was further improved in \cite{balas2004pivot}, allowing variables shifting up and down. Another similar pivot heuristic is the one proposed by Eckstein in \cite{eckstein2007pivot} for 0-1 Integer programs that exploits pivot operations, cuts, and diving heuristics. Berthold \cite{berthold2006primal} proposed two enhanced versions of Simple Rounding, which allow non LP-feasible updates: Rounding and Shifting. In both the algorithms variables are rounded according to a specific scalar score, built with the aim of minimizing the constraints violations and, eventually, retrieving the LP-feasibility following infeasible updates. However, while the first one performs only rounding of integer variables, Shifting adjusts continuous variables as well, shifting their value to specific quantities. Another rounding-based method is RENS, proposed by Berthold in \cite{berthold2007rens} and \cite{berthold2014rens} and based on the solution of smaller MIP sub-problems.

Another important class of start heuristics is the LP-based heuristics, which try to achieve an integer feasible solution by solving a sequence of LP problems. This main idea was introduced by Fischetti and Lodi in \cite{fischetti2005feasibility,fischetti2009feasibility,bertacco2005feasibility} with the Feasibility Pump, further improved in \cite{boland2014boosting,baena2010using} and extended even to nonlinear mixed-integer problems, both convex in \cite{bonami2009feasibility} and non-convex in \cite{d2010experiments}.

The optimization research community provided us with an enormous amount of completely new start heuristics, such as Octane \cite{balas2001octane}, mainly based on geometrical ideas, and an adaptation of Tabu-search \cite{glover1997tabu}, exploited to solve a binary problem by Lokketangen in \cite{lokketangen1998solving}. However, an exhaustive analysis of all the available MIP heuristics is far beyond our work, which is instead mainly focused on developing a new mixed start heuristic; for such an analysis the reader should refer to \cite{berthold2006primal, hanafi2017mathematical}.

Our research is focused on start heuristics and, in particular, on the development of a novel heuristic, PUSH, to find a first feasible point without any further assumption on the problem structure. We modify the Feasibility Pump, introducing randomness to prevent the algorithm from cycling on a finite number of solutions or getting stuck in the same point. Furthermore, we suggest a new rounding rule to update the variables between subsequent LP-optimizations. Instead of rounding to the closest integer, which can easily lead to increase violations, we decide to take into account the weight of each variable in the constraints, performing the Shifting algorithm. This allows either to round integer variables such that the violations are minimized, or even to retrieve the LP-feasibility. Our computational experience on MIPLIB 2017 instances show that, following this approach, we obtain success rate values that are significant higher than with other rounding heuristics. Furthermore, despite an affordable increase in terms of computational time, PUSH generally returns better solutions than only Feasibility Pump.

This paper is organized as follows. In Section \ref{sec:preliminaries_and_notation} we formalize our problem and introduce some notation that will be used further to describe the algorithms. In Section \ref{sec:feasibility_pump}, after describing the rounding heuristics we have used, we discuss the Feasibility Pump, and we explain, how we enhanced this algorithm to avoid the occurrence of cycles. In Section \ref{sec:push} we describe the proposed start heuristic PUSH. Eventually, in Section \ref{sec:computational_experience}, we discuss our computational experience on 800 MIPLIB 2017 instances, and in \ref{sec:Conclusions} we draw our conclusions.


\section{Preliminaries and notation}\label{sec:preliminaries_and_notation}
We address the generic mixed-integer linear problem in the following form:

\begin{equation}
    \label{eq:P}
    P : \left\{ 
\begin{array}{rl}
    \ds\min_{x \in \setR^n} & c^\intercal x \\
     \text{s.t.}& A x \le b\\
     & D x = f\\
     & x_i \in \setZ, \ i \in I
\end{array}
\right.
\end{equation}

We call  $n$  the number of variables,  $n_I$  the number of integer variables (including binaries) and
 $n_C$ the number of continuous variables. The set of indices of the integer variables is $I$, whereas $C$ is the one for the continuous ones. All the coefficients are supposed to be real values, meaning
 $c \in \setR^n$,
  $A\in\setR^{m\times n}$, $b \in \setR^m$, $D\in\setR^{p\times n}$ and $f \in \setR^p$. We will call
     $l_i, u_i \in \setR$   the lower and the upper bounds respectively for the variables $x_i$, $i= 1, \dots, n$.

By dropping out the integer requirements for the variables, we obtain the linear relaxation of $P$, namely
\begin{equation}
    \label{eq:LP}
LP : \left\{ 
\begin{array}{rl}
    \ds\min_{x \in \setR^n} & c^\intercal x \\
     \text{s.t.}& A x \le b\\
     & D x = f
\end{array}
\right.
\end{equation}

{\parindent 0cm
We will also use the rounding  operator  $\lceil a \rfloor : \setR \mapsto \setZ$, defined as $\lceil a \rfloor = \left\{ \begin{array}{ll}
\lfloor a \rfloor &     \text{if }   a - \lfloor a \rfloor \le 0.5         \\
     \lceil a \rceil&   \text{oth.}
\end{array} \right.$
}

For the purpose of evaluating PUSH against other state-of-the-art start heuristics we employ two performance metrics: the CPU time required to obtain a feasible solution or to return a failure outcome, and the optimality gap, computed as follows:
\begin{equation}
    \label{eq:opt_gap}
    \gamma (x) = \frac{c^Tx - c^Tx^*}{c^Tx^*}
\end{equation}
Where $x^*$ is the optimal integer solution or, in case the instance is not solved to optimality yet, the best known solution reported in MIPLIB 2017 documentation \footnote{\url{https://miplib.zib.de/}}.

\section{Feasibility Pump and Rounding techniques}
\subsection{Feasibility Pump}
\label{sec:feasibility_pump}
The Feasibility Pump (FP) is a primal heuristic originally proposed in \cite{fischetti2005feasibility} to find a feasible solution for Mixed Binary Problems, and then extended to MILP in \cite{bertacco2005feasibility} and \cite{achterberg2007improving}.
Given a MIP problem $P$ as in \eqref{eq:P} and its linear relaxation $LP$ as in \eqref{eq:LP}, the polyhedron $Q = \{ x \in \setR^n : Ax \le b, Dx = f\}$ is the feasible set of $LP$.
The Feasibility Pump  proceeds by alternating steps of projection and rounding, starting from the optimal solution  for $LP$, i.e. $x^0 \in \arg\min LP$. The idea is to generate two sequences, one dealing with integer feasibility and the other with $LP$-feasibility. The first sequence   $\{y^t\}$ is obtained by rounding to the nearest integer the values of the points derived step by step from the other sequence, i.e. $y^t = \left\{\begin{array}{cr}
     \lceil x^{t-1}_i\rfloor& \text{ if } i \in I \\
     x^{t-1}_i& \text{oth.}
\end{array} \right.$. 
The second sequence $\{\widetilde{x}^t\}$ is composed by points  belonging to  $Q$, obtained by solving the $L1-$norm projection of the points in the first sequence,i.e. $x^t \in \arg\min \{ ||x - y^t||_1 : x \in Q  \}$. In this way the algorithm eventually terminates when $x^t = y^t$ and the optimal $L1$-norm value is zero. Form the computational perspective, the $L1$-norm function can be linearized, obtaining the following linear program. The outline of Feasibility Pump is reported in algorithm \ref{algo:feasibility_pump}.
\begin{equation}
    \label{eq:L1LP}
    L1LP(y) = \left\{ 
\begin{array}{rl}
    \ds\min_{x \in \setR^n} & \sum_{i \in I} z_i\\
     \text{s.t.}& A x \le b\\
     & D x = f
     \\ & z_i \ge  y_i - x_i, \ i \in I
     \\ & z_i \ge x_i -  y_i , \ i \in I
\end{array}
\right.
\end{equation}

\begin{algorithm}\caption{Basic Feasibility Pump}\label{algo:feasibility_pump}
\DontPrintSemicolon
\SetAlgoLined
Initialization: \textit{maxiter}\;

Let $\widetilde{x}$ be a solution for LP \;

\For{ $t=1, \dots, maxiter$ }{
Set $y : y_i = \left\{
\begin{array}{cc}
    \lceil \widetilde{x}_i \rfloor & \text{if } i \in I \\
     \widetilde{x}_i &  \text{oth.}
\end{array} \right.$\;
\eIf{$L1LP(y) = 0$ }{
Return y \;
}{
Retrieve a solution $(x, z)$ from $L1LP(y)$, i.e. $(x, z) \in \arg L1LP(y)$\;}
}
\end{algorithm}

However, the original version of Feasibility Pump is not proved to converge to a feasible solution under any assumption and the algorithm suffer from two structural drawbacks: poor objective function quality and cycles.
The first issue is a consequence of the fact that, when solving the \eqref{eq:L1LP} at iteration $k$, we do not take into account the initial objective function $c^Tx_k$, i.e., we do not impose any constraint on the objective function decrease, meaning $c^T x_{k+1}$ could easily be further worse than $c^Tx_k$. To partially overcome this issue we introduce the parameter $\alpha$, which is supposed to give importance to the objective function when solving \eqref{eq:L1LP}, whose objective is modified as follows:
$$ f_{L1LP}(z,\alpha;x) =  (1-\alpha)\sum_{i \in I}z_i +\alpha\frac{\sqrt{|I|}}{||c||}c^Tx $$
The closer $\alpha$ to the unity, the higher will be the importance given to the original objective of \eqref{eq:P} and the less the priority given to the minimization of the $L1-$norm.

The second issue, the possibility of ending up in cycles, is even more dangerous, as it directly affects the success rate of the algorithm. Two types of cycles can occur:
\begin{itemize}
    \item One length cycle that happens when $\widetilde{x}^{t} = \widetilde{x}^{t-1}$
    \item $t-$length cycle that happens when $\widetilde{x}^t=\widetilde{x}^{t-k}$
\end{itemize}
For each one of these two types of cycles we implemented a different rule to tackle it. 
For one length cycle we select the $T$ most fractional variables in $\widetilde{x}^{t}$ and for each one we round it to the opposite side, e.g., if $\widetilde{x}^{t}_j = 0.3$ and $j\in T$ we round it to 1 instead to 0.
Concerning $t-$length cycles, for each fractional variable we sample a real number $\rho$ uniformly distributed in $(-0.3;0.7)$ and, if it holds that $|\widetilde{x}_i-\lceil \widetilde{x}_i\rfloor+\max\{0;\rho\}| > 0.5$, we round $\widetilde{x}_i$ to the opposite direction. 
Despite these two rules still do not guarantee the convergence of the algorithm, we nevertheless experiment in many instances that, when adding these two conditions, FP was able to find a solution that could not find in the basic version.

Eventually, we also have to notice that Feasibility Pump shares an efficiency issue with all $LP-$based heuristics; a linear problem of the same dimension than the original MIP problem is solved at every iteration. When dramatically increasing the instance size, this can result in unaffordable computational effort.

\subsection{Simple Rounding and Rounding}\label{sec:simple_rounding}
Rounding methods are a class of start heuristics based on the attempt to round fractional variables while preserving or retrieving LP-feasibility. They can be very efficient in terms of computational time, as they require only basic mathematical operations, but they perform poorly in most medium or large-sized instances. In particular, as already remarked by Berthold in \cite{berthold2006primal}, they suffer from two main drawbacks. First, the rounded variable is usually chosen according to some trivial criteria, such as the number of locks, affecting the algorithm's capability of maintaining LP-feasibility over a large number of iterations. Second, most of the selection criteria are not sensitive to the loss of LP-feasibility, i.e., they do not provide special rules to rapidly retrieve the feasibility after it has been lost for the first time. The heuristics we discuss further have already been implemented in SCIP, see \cite{scip}. 

One of the first rounding techniques is Simple Rounding, see
\cite{berthold2006primal}-\cite{achterberg2007constraint}, developed to tackle problems without equalities, i.e. $p=0$. Starting from an LP-feasible point, the method looks for fractional variables, called trivially roundable,  that can be rounded up or down without losing the feasibility. If there are no such variables, the algorithm stops and returns a failure. The Simple Rounding outline is shown  in Algorithm \ref{ref:SR_algo}.
\begin{algorithm}[h!]
\caption{Simple Rounding}\label{ref:SR_algo}
\DontPrintSemicolon
\SetAlgoLined
Compute $x$ LP-feasible \;
Compute $\bar I = \{j \in I:\ x_j \ \text{is fractional}\}$ \;
 \For{$j \in \bar I$}{
 Set $\bar x = x$ \;
 \If{$x_j$ is trivially down-roundable}{
 $\bar x_j = \lfloor x_j \rfloor $}
 \If{$x_j$ is trivially up-roundable}{
 $\bar x_j = \lceil x_j \rceil $}
 $x = \bar x$
 }
\end{algorithm}

A more complex procedure is the one called Rounding, see \cite{berthold2006primal} and \cite{achterberg2012rounding}, which improves the Simple Rounding scheme by allowing infeasible updates. When an infeasible update occurs, the algorithm does not stop but tries to retrieve feasibility following the rule of locks minimization. A score is associated with each variable, which decreases when the number of locks increases. The variable with the lowest score is chosen and rounded. The underlying idea is to round up variables with a large number of down-locks and to round down variables with a large number of up-locks.

The Rounding is a generalization of Simple Rounding, which performs the same operation, but also allows infeasible updates. Until the input point at each iteration stays feasible, Rounding follows the same rule of Simple Rounding, selecting a variable which is trivially roundable. However, when an infeasible update occurs, the algorithm does not stop, but tries to retrieve feasibility following the rule of locks minimization. In this case each variable is given a score, which decreases when the number of locks increases. The variable with the lowest score is chosen and rounded. The underlying idea is to round up variables with a large number of down-locks and to round down variables with a large number of up-locks, relying on the fact that rounding down a variable with a lot of up-locks (or, symmetrically, rounding up a variable with a lot of down-locks) means improving constraints violation. However, this is not always true. In fact, our computational experience shows that Rounding alone rarely succeeds in retrieving LP-feasibility. This result is not surprising, as there is no mathematical law directly linking the number of locks with the entity of constraints violation. Furthermore, the number of locks seems to be a quite poor representation of the influence we can have on violations when rounding a variable; for the purpose of weighing the influence on constraints, one should also consider available slacks as well as the normalized magnitude of the locks. The outline of Rounding is reported in Algorithm \ref{ref:R_algo}.

\begin{algorithm}[h!]
\caption{Rounding}\label{ref:R_algo}
\DontPrintSemicolon
\SetAlgoLined
Compute $x$ LP-feasible \;
Compute $\bar I = \{j \in I:\ x_j \ \text{is fractional}\}$ \;
 \While{$|\bar I| > 0$}{
 Set $\bar x = x$ \;
 Check $\bar x$ for LP-feasibility \;
 \If{$\bar x$ is LP-feasible}{
 Set $\xi_{max} = -1$, $j_{min} = \infty$, $\sigma = 0$ \;
  \For{$j \in \bar I$}{
  Compute $\xi_j$ = Number of up-locks of $x_j$ \;
 \If{$\xi_j > \xi_{max}$}{
 $\xi_{max} = \xi_j$, $\sigma = -1$, $j_{min} = j$}
 Compute $\xi_j$ = Number of down-locks of $x_j$ \;
  \If{$\xi_j > \xi_{max}$}{
 $\xi_{max} = \xi_j$, $\sigma = 1$, $j_{min} = j$}
 }
 }
 \Else{
 Choose a violated constraint $i:\ A_i^T \bar x > b_i$ \;
 Set $\xi_{min} = \infty$, $j_{min} = \infty$, $\sigma = 0$ \;
 \For{$j \in \bar I$}{
 \If{$a_{ij} > 0$}{
 Compute $\xi_j$ = Number of down-locks of $x_j$ \;
 \If{$\xi_j < \xi_{min}$}{
 $\xi_{min} = \xi_j$, $\sigma = -1$, $j_{min} = j$}
 }
  \If{$a_{ij} < 0$}{
 Compute $\xi_j$ = Number of up-locks of $x_j$ \;
 \If{$\xi_j < \xi_{min}$}{
 $\xi_{min} = \xi_j$, $\sigma = 1$, $j_{min} = j$}
 }
 }
 }
 \If{$\sigma = 1$}{ $\bar x_j = \lceil x_j \rceil $}
 \Else{$\bar x_j = \lfloor x_j \rfloor $}
 Set $x = \bar x$
 }
\end{algorithm}

\subsection{Shifting}\label{sec:shifting}
In \cite{berthold2006primal}, Berthold introduced the Shifting algorithm, a rounding algorithm that can also update integer or continuous variables multiple times to retrieve LP-feasibility.

Shifting performs Simple Rounding iterations, whenever the current point is LP-feasible and looks for a particular variable to update otherwise. While during a Rounding iteration the variable to update is always chosen in the set $\bar I$ of the fractional variables in $I$, Shifting also allows to update both integer variables in $I \setminus \bar I$ and continuous variables in order to retrieve LP-feasibility.

Once a violated constraint is selected, each variable, including the ones is in $C$, is given a score depending on the number of locks. The scoring policy facilitates fractional variables in $ I$, but also the already allocated integer ones or the continuous ones can be selected. If a previously rounded variable is selected, then it is decreased or increased by one according to the update sign. If a continuous variable is selected, it is shifted to a certain quantity to decrease the current constraint violation. Shifting a variable that can take any value within the feasible region is often an efficient way to rapidly decrease the constraint violation without modifying any integer element. Among the several existing shifting rules, a famous one is the optimal shifting concerning the selected constraint. Provided the algorithm decides to reduce the violation of constraint $i$ through variable $x_{\bar j} \in C$, the optimal shifted value is
$ x^*_j = x^{old}_{\bar j} - \frac{V_i}{a_{i\bar j}}
$.
Where $V_i = \sum_{h=1}^n a_{ij} x^{old}_j - b_i $, is the violation of the $i$-th constraint. Thus, according to this rule, the shifted variable is $x^*_j$ if $x^*_j \in [l_j,u_j]$, otherwise $x_j$ is set to closest feasible value to $x^*_j$.
We have also added an anti-cycle rule to prevent the algorithm from being stuck when a variable, that is already set to its bound, is selected. When such a case occurs, the variable is stored in a list of forbidden variables for the following 50 iterations or until it is updated in the opposite direction with respect to its bound.
The backbone of the Shifting is reported in algorithm \ref{ref:SH_algo}.
\begin{algorithm}[h!] \caption{Shifting}\label{ref:SH_algo}
\DontPrintSemicolon
\SetAlgoLined
Compute $x$ LP-feasible, $\bar I = \{j \in I:\ x_j \ \text{is fractional}\}$ \;
 \While{$|\bar I| > 0$ or $x$ is not LP-feasible}{
 Set $\bar x = x$ \;
 \uIf{$\bar x$ is LP-feasible}{
 Do as Rounding \;
 } \Else{
 Choose a violated constraint $i:\ A_i^T \bar x > b_i$ \;
 Set $\zeta_{min} = \infty$, $j_{min} = \infty$, $\sigma = 0$ \;
 \For{$j:\ a_{ij} \neq 0$}{
 \uIf{$a_{ij} < 0$}{
 Compute $\xi_j$ = Number of down-locks of $x_j$, and set $\zeta_j = -1 + \frac{1}{\xi_j + 1}$ \;
 \If{$\zeta_j < \zeta_{min}$}{
 $\zeta_{min} = \zeta_j$, $\sigma = -1$, $j_{min} = j$}
 }
 \If{$a_{ij} > 0$}{
 Compute $\xi_j$ = Number of up-locks of $x_j$,  and set $\zeta_j = -1 + \frac{1}{\xi_j + 1}$ \;
 \If{$\zeta_j < \zeta_{min}$}{
 $\zeta_{min} = \zeta_j$, $\sigma = 1$, $j_{min} = j$}
 }
 }
 \uIf{$j \in \bar I$}{
 Update according to Rounding rule}
 \Else{
 $\bar x_j = x_j +\sigma$\;
 }
 }
 Set $x = \bar x$
 }
\end{algorithm}

Our computational experience shows that Shifting generally achieves higher success value rate than the other rounding heuristics. Nonetheless, this seems not to have any kind of influence on the solution quality, whereas the complex rounding and shifting operations can take much longer CPU time than in Rounding and Simple Rounding, leading to an efficiency issue.

\section{Proposed method: PUSH}\label{sec:push}
The approach we present here is called PUSH, directly derived from the feasibility PUmp (section \ref{sec:feasibility_pump}) and Shifting (section \ref{sec:shifting}) approaches. The basic idea behind the algorithm is to replace the rounding stage in the feasibility pump with several rounds of rounding. We aim to help some drawbacks of both procedures by combining them, enriching the decision steps with additional measures.

In algorithm \ref{algo:push}, we report the backbone of the method. The   checks for feasibility mostly rely on the Feasibility Pump approach, whereas for rounding and perturbing we developed two algorithms based on the rounding techniques described before.

\begin{algorithm}[h!]
\caption{PUSH}\label{algo:push}
\DontPrintSemicolon
\footnotesize

\SetAlgoLined
Parameters: \textit{maxiter}, \textit{rounding\_threshold}, \textit{shrinking\_trials}, \textit{random\_sensitivity}\;
Initialization:  $\mathcal{L} = \emptyset$\;

Let $x$ be a solution for LP \;

\For{ $t=1, \dots, maxiter$ }{
Set $y = \texttt{push\_rounding}(x,  rounding\_threshold)$\;
\If{ $y \in \cal L$}{
$y = \texttt{random\_perturbation}(y, z, \bar{I},  \mathcal{L}, \textit{random\_sensitivity})$\;
}
$\mathcal{L} \rightarrow \mathcal{L}\cup \{ y \}$\;

\eIf{$L1LP(y) = 0$ }{
Return y \;
}{
Retrieve a solution $(x, z)$ from $L1LP(y)$, i.e. $(x, z) \in \arg L1LP(y)$\;}
}

\end{algorithm}

\subsection{Random perturbation}
In most of the versions of the Feasibility Pump, random perturbations are introduced to avoid unwanted cycling conditions. In particular, it is possible for the method to generate a sequence of couples projection-rounding, that cyclically returns to themselves. As described in section \ref{sec:feasibility_pump}, a typical way is to check if the rounded point has been already exploited, and if it had then a random perturbation is invoked until a new unseen point is found. Inside the PUSH algorithm, we introduce a special way of perturbing the point, that keeps into account some measure of the distance from the point itself, reducing the space of uncertainty and favoring closer solutions.
 
 Taking the L1LP formulation from equation \eqref{eq:L1LP}, we may notice that the $z$-variables induce implicit bounds on  $x$, i.e. $y_i - z_i  \le  x_i   \le y_i + z_i$, $ i \in I$. Despite being significant in reducing the size of the sampling space, this interval may be too strict, especially in the case of binary variables. To this aim we take the lower and upper rounding on both sides, cutting them according to the bounds. Therefore, the single variables are sampled according to a discrete uniform distribution $\mathcal{U} \left\{ \max\{l_i, \lfloor y_i - z_i \rfloor, \min\{u_i, \lceil y_i - z_i \rceil\}\right\}$.
 
Moreover, perturbing all the integer variables may result in unwanted or bad solutions, even if taken according to the rule above. For this reason, we enforce the perturbation on the rounded variables in $\bar I$ which are the ones that were fractional, and we keep a reduced variability for all the others. The motivation is not to deteriorate already integer variables. We introduce a hyper-parameter called \textit{random\_sensitivity} which represents the probability of perturbing a variable whose index is in $\bar I$, whereas the probability of perturbing any other integer variable is ten times smaller, i.e. $\frac{\textit{random\_sensitivity}}{10}$.

The method is reported in algorithm \ref{algo:random_perturbation}.

\begin{algorithm}[h!]
\caption{\texttt{random\_perturbation}}\label{algo:random_perturbation}
\DontPrintSemicolon
\footnotesize

\SetAlgoLined
Input: $y, z, \bar{I},  \mathcal{L}, \textit{random\_sensitivity}$\;

\While{True}{
Compute $L : L_i =\left\{ \begin{array}{ll}
     \max\{l_i, \lfloor y_i - z_i \rfloor\}& \text{ if } i \in 
     \bar I\\
     \max\{ l_i, y_i \}&  \text{ oth.}
\end{array} \right.$, ${U} : {U}_i =\left\{ \begin{array}{ll}
     \min\{u_i,\lceil y_i + z_i \rceil \}& \text{ if } i \in 
     \bar I\\
     \min\{ u_i, y_i \}&  \text{ oth.}
\end{array} \right.$\; 
Compute $\hat y : \hat y_i = \left\{ 
\begin{array}{ll}
    \sim {\cal U} \left\{ L_i, {U}_i \right\} & \text{with probability } \textit{random\_sensitivity} \textit{ and }  \text{ if } i \in 
     \bar I  \\ 
      \sim {\cal U} \left\{ L_i, U_i \right\} & \text{with probability } \frac{ \textit{random\_sensitivity}}{ 10 } \textit{ and } \text{ if } i \not\in
     \bar I  \\ 
     y_i & \text{ oth.}
\end{array}
\right.$\;
\eIf{$\hat y \in \cal L$}{
Continue \;}{
Return $\hat y$\;
}
}
\end{algorithm}

\subsection{PUSH rounding}
The idea behind the rounding procedure in the PUSH algorithm is to take advantage of known rounding methods like Simple rounding (section \ref{sec:simple_rounding}) and Shifting (section \ref{sec:shifting}) to overcome the straightforward mechanism adopted in the Feasibility Pump. Our version of the rounding heuristic focuses on the integer variables, delegating to the projection step the assignment for the continuous ones.

As reported in algorithm \ref{algo:push_rounding}, the procedure is composed of an outer loop that counts the number of variables to round, according to the hyper-parameter \textit{rounding\_threshold}, which we will further indicate as \texttt{rt}.  For every step of the loop, similar to the Shifting algorithm (section \ref{sec:shifting}), our rounding heuristic searches for violated constraints. To improve efficiency, we decided to check in the pool of constraints randomly, giving the greater probability to the constraints that were violated in the previous iterations. When the method finds a violated constraint, it computes the set of candidates, i.e. the set of fractional variables with a non-zero coefficient in the constraint. The intuition is that hard constraints will be violated more frequently than others. The hyper-parameter \texttt{rt} plays a crucial role in balancing the trade-off between computational efficiency and success rate. In fact, the higher the rounding threshold, the higher is also the probability to perform a smart rounding of all the variables and to end up in a feasible point. Despite one could argue that fixing an indiscriminately high rounding threshold maximizes the success rate, we need to consider than if we fix $\texttt{rt} > 1$, we perform more than one rounding per each variable, meaning that we are basically apply the Shifting algorithm between two Feasibility Pump iteration. Such an approach would be affected by the same efficiency issues encountered in Shifting.

\begin{algorithm}[h!]
\caption{\texttt{push\_rounding}}\label{algo:push_rounding}
\DontPrintSemicolon
\footnotesize
\SetAlgoLined
Input: x,  \textit{rounding\_threshold}\;
Initialization: \textit{violated\_constraint} = None\;
Compute $\bar I = \{j \in I:\ x_j \ \text{is fractional}\}$ \;
\For{$t=1, \dots, \lceil n \cdot \textit{rounding\_threshold}\rfloor$}{
    \If{$\bar I = \emptyset$}{
         Return x\;}
    \For{$\textit{constraint} \in \textit{constraints}$, chosen randomly}{
         Find involved variables $\textit{vars} = \{i : \textit{constraint}_i \neq 0, i \in I\}$\;
         \If{$\textit{constraint}$ is violated}{
            Compute $\textit{candidates} = \bar I \cap \textit{vars}$\;
            Assign $\textit{violated\_constraint} = \textit{constraint}$\;
            Break\;
         }
    }
    
    \uIf{\textit{violated\_constraint} is not None and $\textit{candidates} \neq \emptyset$}{
        Choose $\tilde i \in \bar I$,  $d \in \{\textit{down, up}\}$ via \texttt{scoring\_selection}$(\textit{candidates}, \textit{constraint})$\;
    }\uElseIf{ \textit{violated\_constraint} is not None and $\textit{candidates} = \emptyset$}{
        Choose $\tilde i \in \bar I$,  $d \in \{\textit{down, up}\}$ via \texttt{simple\_rounding\_selection}$(\bar I )$\;
    }\Else{
        Choose $\tilde i \in \bar I$, $d \in \{\textit{down, up}\}$ via \texttt{feasible\_selection}$(x, \bar I )$\;
    }
    
    Update point $x \rightarrow \tilde x : \tilde x_i = \left\{ \begin{array}{cl}
        \lfloor x_i \rfloor & \text{ if } i = \tilde i \text{ and } d = down \\
        \lceil x_i \rceil & \text{ if } i = \tilde i \text{ and } d = up\\
        x_i & \text{oth.}
    \end{array}\right.$\; 
    Update $\bar I \rightarrow \bar I \setminus \{\tilde i \}$\;
}
\end{algorithm}

After the search for a violated constraint is concluded, we adopt different strategies to decide which variable to round, according to the specific situation we ended up with. In particular, we may encounter three situations:
\begin{itemize}
    \item[$\cdot$] the current point is infeasible and there are candidates
    \item[$\cdot$] the current point is infeasible but there are no candidates
    \item[$\cdot$] the current point is feasible.
\end{itemize}

The \texttt{scoring\_selection} method, presented in algorithm \ref{algo:scoring_selection}, tackle the first case. For every fractional variable among the candidates, a rounding direction is chosen according to the value of the coefficient in the constraint. Differently from the basic Shifting version, here we consider directly equality constraints by taking into account their slacks. In this way, if a slack and the corresponding coefficient are both positive, then it would be reasonable to round down the variable to reduce the distance from zero. Of course, if the slack is negative for an equality constraint, we adopt the opposite strategy.  

After that, we use the rounding direction to generate an appropriate score. We decided to use a multiplicative score that takes into account the number of locks and the relative magnitude of the variable, represented by the average and normalized coefficient that the variable takes concerning the others. Called $\#locks$ the number of locks and $magnitude$ the magnitude before described, the final score is 
\begin{equation}
    \label{eq:score}
    score = \#locks \cdot e^{magnitude}
\end{equation}
Once all the scores are collected, the one with the highest value is selected.

\begin{algorithm}[h!]
\caption{\texttt{scoring\_selection}}\label{algo:scoring_selection}
\DontPrintSemicolon
\footnotesize

\SetAlgoLined
Input: \textit{cadidates, constraint}\;
Initialization: $\mathcal{S} = \emptyset$, $\tilde i$, $d$\;

\For{$ j \in candidates$}{
        \uIf{\textit{constraint} is an inequality and $\textit{sign}\left( \textit{constraint}_j \right) > 0$}{
            Set $d_j = \textit{down}$\;
        }\uElseIf{\textit{constraint} is an equality and $\textit{sign}\left( \textit{constraint}_j \right) \cdot \textit{sign}(constraint.slack) > 0$}{
            Set $d_j = \textit{down}$\;
        }\Else{
            Set $d_j = \textit{up}$\;
        }
        
        \eIf{$d_j = \textit{down}$}{
            Compute down locks $\textit{down}_j$ and  magnitude $\text{mag\_down}_j$ for $x_j$\;
            Update $\mathcal{S} \rightarrow \mathcal{S} \cup \{\textit{down}_j \cdot e^{\textit{mag\_down}_j} \}$\;
        }{
            Compute up locks $\textit{up}_j$ and  magnitude $\text{mag\_up}_j$ for $x_j$\;
            Update $\mathcal{S} \rightarrow \mathcal{S} \cup \{\textit{up}_j \cdot e^{\textit{mag\_up}_j} \}$\;
        }

}

Assign $\tilde i \rightarrow \arg\max \cal S$, $d \rightarrow d_{\tilde i}$\;
Return $\tilde i$, $d$\;
\end{algorithm}

It may happen that the set of fractional variables does not contain indices connected to violated points. In this case, the set of candidates is empty and the rounding procedure invokes the \texttt{simple\_rounding\_selection} method, which is based on the Simple Rounding technique, see sec.\ref{sec:simple_rounding}. Every fractional variable is inspected, counting the number of up and down locks. If a variable is trivially roundable, then the algorithm returns the index and the rounding direction accordingly, otherwise, it collects scores in the form of \eqref{eq:score}. Finally, the variable with the highest score is returned, and the direction is chosen depending on whether the score was generated by counting up or down values.

\begin{algorithm}[h!]
\caption{\texttt{simple\_rounding\_selection}}\label{algo:simple_rounding_selection}
\DontPrintSemicolon
\footnotesize

\SetAlgoLined
Input: $\bar I$\;
Initialization: $\textit{DOWN} = \textit{UP} = \emptyset$, $\tilde i$, $d$\;

\For{$j \in \bar I$}{
    Compute down $\textit{down}_j$ and up $\textit{up}_j$ locks  for $x_j$ \;
    \uIf{$\textit{down}_j = 0$ and $\textit{up}_j > 0$}{
        Assign $\tilde i \rightarrow j$, $d \rightarrow \textit{up}$\;
        Return $\tilde i$, $d$\;
    }\ElseIf{$\textit{up}_j = 0$ and $\textit{down}_j > 0$}{
        Assign $\tilde i \rightarrow j$, $d \rightarrow \textit{down}$\;
        Return $\tilde i$, $d$\;
    }
    Compute normalized magnitudes for $x_j$, resp. $\textit{mag\_down}_j$ and $\text{mag\_up}_j$\;
    Update $\textit{DOWN} \rightarrow \textit{DOWN} \cup \{\textit{down}_j \cdot e^{\textit{mag\_down}_j} \}$\;
    Update $\textit{UP} \rightarrow \textit{UP} \cup \{\textit{up}_j \cdot e^{\textit{mag\_up}_j} \}$\;
}

Compute $j_{down} = \arg\max \textit{DOWN}$, $\textit{max\_down\_score} = \max \textit{DOWN}$\;
Compute $j_{up} = \arg\max \textit{UP}$, $\textit{max\_up\_score} = \max \textit{UP}$\;
\eIf{$\textit{max\_up\_score} > \textit{max\_down\_score}$}{
    Assign $\tilde i \rightarrow j_{up}$, $d \rightarrow \textit{up}$\;
}{
    Assign $\tilde i \rightarrow j_{down}$, $d \rightarrow \textit{down}$\;
}

Return $\tilde i$, $d$\;
\end{algorithm}

The last and rare case is the one where the point is feasible, but some fractional variables are still present. The rounding algorithm chooses the variable and the direction by calling the \texttt{feasible\_selection} method. Similarly to before, our final goal is to compute a score and take the best out, but instead of taking the contribution concerning violated constraints, we consider the residuals. The reason behind this is to facilitate the selection of those variables influencing not too tight constraints. After collecting all the absolute values of the residuals, and dividing them between the one where the variable has a positive coefficient or not, the algorithm computes a new multiplicative score in the form
$$
    score = \#locks \cdot e^{\frac{magnitude}{\textit{avg\_res}}}
$$
where \textit{avg\_res} is the average value of the residuals.

\begin{algorithm}[h!] \caption{\texttt{feasible\_selection}}\label{algo:feasible_selection}
\DontPrintSemicolon
\footnotesize

\SetAlgoLined
Input: $x$, $\bar I$\;
Initialization: $\mathcal{S} = \emptyset$\;

\For{$j \in \bar I$, chosen randomly}{
     Set $\textit{residuals\_down}= \textit{residuals\_down}= \emptyset$\;

     Compute the residuals of the constraints where the  $j-$th coefficient is positive, then store their absolute value in \textit{residuals\_down}\;
     Compute the residuals of the constraints where the  $j-$th coefficient is negative,  then store their absolute value in \textit{residuals\_up}\;
     
     
     

     Compute down locks and down normalized magnitude for $x_j$, resp. $\textit{down}_j$ and $\text{mag\_down}_j$\;
    Compute up locks and up normalized magnitude for the variable $x_j$, resp. $\textit{up}_j$ and $\text{mag\_up}_j$\;
                
    Compute the average residuals, rep. \textit{avg\_down\_res} and \textit{avg\_up\_res}\;

    Set $\textit{score\_down}=  \textit{down}_j \cdot \exp{\frac{\textit{mag\_down}_j}{\textit{avg\_down\_res}}}$\;
    Set $\textit{score\_up}=  \textit{up}_j \cdot \exp{\frac{\textit{mag\_up}_j}{\textit{avg\_up\_res}}}$\;
    
    \eIf{$\textit{score\_down}>\textit{score\_up}$}{
        Set $d_j =down$, $\mathcal{S} \rightarrow \mathcal{S} \cup \left\{ \textit{score\_down}\right\}$\;
    }{
        Set $d_j =up$, $\mathcal{S} \rightarrow \mathcal{S} \cup \left\{ \textit{score\_up}\right\}$\;
    }

}

Assign $\tilde i \rightarrow \arg\max \cal S$, $d \rightarrow d_{\tilde i}$\;
Return $\tilde i$, $d$\;
\end{algorithm}

\section{Computational Experience}\label{sec:computational_experience}
We conducted our tests on an Intel(R) Core(TM) i5-10210U CPU @ 1.60GHz 2.11 GHz 16 GB. The code is written in Python v3.8 and it uses Gurobi  v9.5.0  \cite{gurobi} to solve the linear programs. We tested PUSH against Simple Rounding, Rounding, Shifting and Feasibility Pump using 800 test instances, that were taken from the open-source MIPLIB 2017 \cite{miplib} library.
For the purpose of performing a fair comparison, we firstly selected a subset of easy instances taken from MIPLIB 2017 \cite{miplib} to tune the rounding threshold $\texttt{rt}$ in PUSH and, after carrying out the tests reported in table \ref{tab:PUSH_setting_rt} to find an optimal value, we set  $\texttt{rt} = 0.6$. 

We further report, in table \ref{tab:Overview}, the overall performances of the five algorithms, where the average gap is computed according to \eqref{eq:opt_gap} excluding a certain number of instances, which are considered outliers as the value of the gap obtained in those instances is not comparable with the average, and only represents added noise. Looking at \eqref{eq:opt_gap}, it is clear that this can happen when the optimal value of the objective function is either $0$ or very close to $0$. Furthermore, the average gap does not include unbounded and infeasible instances, for which we cannot define a gap. Nevertheless, for the sake of completeness, the results with outliers included can be found in the Appendix.

Results reported in table \ref{tab:Overview} show that the rounding heuristics are largely outperformed by Feasibility Pump and PUSH. Despite PUSH does not seem to have a significant advantage over Feasibility Pump in terms of success rate, the average gap appears to be almost 10 \% lower. 

We also wanted to verify, whether this behaviour is stable when varying the problem complexity, and to assess PUSH efficiency in terms of average computational time. For this purpose we split the 800 instances into four quartiles, ordered by number of non-zeros, i.e., the number of coefficients $a_{ij} \ne 0$ of matrix $\left[A \ \ \ D\right]^T$. For every quartile we report in table \ref{tab:PUSH_quartili} the average performances and the average problem dimension in terms of variables, constraints and non-zeros. We observe that PUSH outperforms Feasibility Pump in terms of average gap in all the four quartiles. Furthermore, PUSH seems not to require much more computational effort than Feasibility Pump and in the second quartile PUSH is even more efficient. This can be explained considering that PUSH needs more time to decide which variables should be rounded but at the same time it generally requires significantly less iterative solution of \eqref{eq:L1LP}.

\begin{table}[h!]
\footnotesize
    \centering
    \caption{PUSH performance on a sample dataset from MIPLIB2017 with different \texttt{rt}. The symbol $-$ indicates  the algorithm was not able to find a feasible solution after 250 iterations.}
    \label{tab:PUSH_setting_rt}
    \begin{tabular}{lrrrrrrrrrr}
    \hline
     Instance&\multicolumn{2}{c}{\texttt{rt} = 0.2 } &  \multicolumn{2}{c}{\texttt{rt} = 0.4 }&  \multicolumn{2}{c}{\texttt{rt} = 0.6 }&  \multicolumn{2}{c}{\texttt{rt} = 0.8 }&  \multicolumn{2}{c}{\texttt{rt} = 1 } \\ \hline

         & Time & Obj & Time & Obj & Time & Obj & Time & Obj & Time & Obj \\ \cline{2-11}
        10teams  & 70.8 & 1009.0 & 5.08 & 996 & 32.6 & 1019.3 & 37.4 & 1020.7 & 46.5 & 102- \\ 
        22433  & 0.5 & 21608 & 0.17 & 21583 & 0.45 & 21559 & 0.17 & 21578.6 & 0.6 & 21603 \\ 
        23588  & 0.3 & 8326.2 & 0.13 & 8343 & 0.6 & 8302.0 & 0.6 & 8304.4 & 0.4 & 8313 \\ 
        50v-10  & 0.3 & 440497 & 0.3 & 579080 & 0.3 & 308675 & 0.4 & 435898 & 0.5 & 541746 \\ 
        acc-tight2  & 8.6 & 0.0 & 6.3 & 0.0 & 5.2 & 0.0 & 9.5 & 0.0 & 8.7 & 0.0 \\ 
        air03  & 6.0 & 379220 & 7.3 & 3.89E05 & 5.6 & 383637 & 9.4 & 426222 & 8.7 & 385214 \\ 
        air05  & 7.7 & 35522 & 4.9 & 34066 & 6.6 & 36134 & 11.2 & 37361 & 11.3 & 34255 \\ 
        arki001  & - & - & - & - & - & - & - & - & - & - \\ 
        comp07-2idx  & 141.5 & 9472.4 & 138.6 & 9819 & 172.9 & 9179.5 & 178.7 & 9412.4 & 174.5 & 8878.4 \\ 
        diameterc-msts  & 1130.4 & 729.0 & 463.1 & 729.0 & 670.7 & 729.0 & 1175.3 & 732.0 & 390.5 & 729.0 \\ 
        eil33-2  & 19.2 & 2137.9 & 12.4 & 2075.3 & 22.8 & 2080.3 & 16.1 & 1784.3 & 35.7 & 1870.3 \\ 
        g503inf  & - & - & - & - & - & - & - & - & - & - \\ 
        gen-ip054  & 0.05 & 1.1E11 & 0.06 & 1.9E11 & 0.0 & 6.1E10 & 0.0 & 3.0E10 & 0.0 & 5.9E10 \\ 
        mcsched  & 0.5 & 442123 & 0.4 & 442318 & 0.4 & 442318 & 0.5 & 442129 & 0.5 & 442533 \\ 
        mod010  & 0.6 & 6864 & 0.6 & 6864 & 0.6 & 6825.6 & 0.7 & 6833.2 & 0.7 & 6839.0 \\ 
        neos-1354092  & - & - & - & - & - & - & - & - & - & - \\ 
        neos-1599274  & 2.5 & 49201 & 1.5 & 55142 & 1.4 & 51731 & 2.7 & 50870 & 1.8 & 48778 \\ 
        neos-2624317-amur  & - & - & - & - & - & - & - & - & - & - \\ 
        neos-3024952-loue  & - & - & - & - & - & - & - & - & - & - \\ 
        peg-solitaire-a3  & - & - & - & - & - & - & - & - & - & - \\ 
        qap10  & 492.1 & 504 & 112.3 & 516 & 67.3 & 485.5 & 188.3 & 504 & 220.1 & 533 \\ 
        r50x360  & 0.6 & 36568 & 0.3 & 27442 & 0.3 & 30352 & 0.4 & 28372 & 0.6 & 30975 \\ 
        rococoC10-001000  & - & - & - & - & - & - & - & - & - & - \\ 
        seymour  & 1.8 & 578.2 & 1.97 & 591 & 1.7 & 556.0 & 2.3 & 590.4 & 2.1 & 571.0 \\ 
        supportcase17  & - & - & - & - & - & - & - & - & - & - \\ 
        xmas10-2  & - & - & - & - & - & - & - & - & - & - \\ \hline
    \end{tabular}
\end{table}

\begin{table}[!ht]
    \centering
    \caption{Performances on 800 MIPLIB 2017 instances}
    \label{tab:Overview}
    \begin{tabular}{|l|l|l|l|l|l|}
    \hline
         Algorithm & Simple Rounding & Rounding & Shifting & Feasibility Pump & PUSH \\ \hline
        $N_{solved}$& 96 & 130 & 156 & 383 & 393 \\
        $N_{outlier}$ & 11 & 16 & 21 & 74 & 70 \\
        Avg Gap & 1.84 & 1.69 & 1.62 & 7.23 & 6.63 \\
        Success & $12\%$ & $16.25\%$ & $19.5\%$ & $47.87\%$ & $49.13\%$ \\ \hline
    \end{tabular}
\end{table}

\begin{table}[h!]
\footnotesize
    \centering
    \caption{Average PUSH and Feasibility Pump performances by quartiles. Infeasible, unbounded instances, as well as outliers are not included}
    \label{tab:PUSH_quartili}
    \begin{tabular}{lrrrrrrrrrr}
    \hline
     Quartile &\multicolumn{3}{c}{Feasibility Pump} &  \multicolumn{3}{c}{PUSH}&  $\bar n_{vars}$& $\bar n_{constrs}$&$\bar n_{nonzeros}$\\ \hline

         & Avg Gap & Avg Time & $N_{solved}$ & Avg Gap & Avg Time & $N_{solved}$ &\\ \cline{2-7}
         1st & 5.71 & 16.99 & 100 & 5.3 & 23.02 & 106 & 1167 & 773 & 4148 \\
        2nd & 7.07 & 23.81 & 76 & 5.7 & 20.88 & 77 & 5825 &  4951 & 24027 \\
        3rd & 6.43 & 36.35 & 79 & 5.62 & 46.59 & 80 & 15060 & 16715&   83296 \\
        4th & 9.98 & 88.42 & 64 & 9.26 & 79.31 & 64 & 41282 &  47696& 274776 \\ \hline 
    \end{tabular}
\end{table}

\section{Conclusions}
\label{sec:Conclusions}
Our computational experience, as well as the test on a huge amount of instances we have carried out, show that we have developed a novel feasible heuristic for Mixed-Integer Programming problems. Despite PUSH does not have any convergence guarantee, our tests show that it generally leads to better solutions than Feasibility Pump and does not require much more computational effort; the higher effort to round the current solution is compensated by a gain in terms of number of Feasibility Pump iterations. Furthermore, PUSH is robust when increasing the instance dimension and the advantage over Feasibility Pump does not seem to be influenced by the instance type or dimension.
Concerning future research perspective, we believe that PUSH can be enhanced, introducing a structural optimized method to choose the violated constraint and the variable to round, possibly using a learning mechanism. Eventually, PUSH could be used jointly with an improvement heuristic to rapidly provide a feasible solution that is not too far from the optimal one.
\clearpage
\bibliographystyle{apalike}
\bibliography{biblio}

\begin{thebibliography}{}

\bibitem[Achterberg, 2007]{achterberg2007constraint}
Achterberg, T. (2007).
\newblock Constraint integer programming.

\bibitem[Achterberg and Berthold, 2007]{achterberg2007improving}
Achterberg, T. and Berthold, T. (2007).
\newblock Improving the feasibility pump.
\newblock {\em Discrete Optimization}, 4(1):77--86.

\bibitem[Achterberg et~al., 2012]{achterberg2012rounding}
Achterberg, T., Berthold, T., and Hendel, G. (2012).
\newblock Rounding and propagation heuristics for mixed integer programming.
\newblock In {\em Operations research proceedings 2011}, pages 71--76.
  Springer.

\bibitem[Baena~Mirabete and Castro~P{\'e}rez, 2010]{baena2010using}
Baena~Mirabete, D. and Castro~P{\'e}rez, J. (2010).
\newblock Using the analytic center in the feasibility pump.

\bibitem[Balas et~al., 2001]{balas2001octane}
Balas, E., Ceria, S., Dawande, M., Margot, F., and Pataki, G. (2001).
\newblock Octane: A new heuristic for pure 0--1 programs.
\newblock {\em Operations Research}, 49(2):207--225.

\bibitem[Balas and Martin, 1986]{balas1986pivot}
Balas, E. and Martin, C. (1986).
\newblock Pivot and shift--a heuristic for mixed integer programming.
\newblock {\em GSIA, Carnegie Mellon University}.

\bibitem[Balas and Martin, 1980]{balas1980pivot}
Balas, E. and Martin, C.~H. (1980).
\newblock Pivot and complement--a heuristic for 0-1 programming.
\newblock {\em Management science}, 26(1):86--96.

\bibitem[Balas et~al., 2004]{balas2004pivot}
Balas, E., Schmieta, S., and Wallace, C. (2004).
\newblock Pivot and shift—a mixed integer programming heuristic.
\newblock {\em Discrete Optimization}, 1(1):3--12.

\bibitem[Bertacco, 2005]{bertacco2005feasibility}
Bertacco, L. (2005).
\newblock A feasibility pump heuristic for general mixed-integer problems livio
  bertacco, matteo fischetti◦, andrea lodi.

\bibitem[Berthold, 2006]{berthold2006primal}
Berthold, T. (2006).
\newblock {\em Primal heuristics for mixed integer programs}.
\newblock PhD thesis, Zuse Institute Berlin (ZIB).

\bibitem[Berthold, 2007]{berthold2007rens}
Berthold, T. (2007).
\newblock Rens-relaxation enforced neighborhood search.

\bibitem[Berthold, 2014]{berthold2014rens}
Berthold, T. (2014).
\newblock Rens.
\newblock {\em Mathematical Programming Computation}, 6(1):33--54.

\bibitem[Bestuzheva et~al., 2021]{scip}
Bestuzheva, K., Besan{\c{c}}on, M., Chen, W.-K., Chmiela, A., Donkiewicz, T.,
  van Doornmalen, J., Eifler, L., Gaul, O., Gamrath, G., Gleixner, A.,
  Gottwald, L., Graczyk, C., Halbig, K., Hoen, A., Hojny, C., van~der Hulst,
  R., Koch, T., L{\"u}bbecke, M., Maher, S.~J., Matter, F., M{\"u}hmer, E.,
  M{\"u}ller, B., Pfetsch, M.~E., Rehfeldt, D., Schlein, S., Schl{\"o}sser, F.,
  Serrano, F., Shinano, Y., Sofranac, B., Turner, M., Vigerske, S.,
  Wegscheider, F., Wellner, P., Weninger, D., and Witzig, J. (2021).
\newblock {The SCIP Optimization Suite 8.0}.
\newblock Technical report, Optimization Online.

\bibitem[Boland et~al., 2014]{boland2014boosting}
Boland, N.~L., Eberhard, A.~C., Engineer, F.~G., Fischetti, M., Savelsbergh,
  M.~W., and Tsoukalas, A. (2014).
\newblock Boosting the feasibility pump.
\newblock {\em Mathematical Programming Computation}, 6(3):255--279.

\bibitem[Bonami et~al., 2009]{bonami2009feasibility}
Bonami, P., Cornu{\'e}jols, G., Lodi, A., and Margot, F. (2009).
\newblock A feasibility pump for mixed integer nonlinear programs.
\newblock {\em Mathematical Programming}, 119(2):331--352.

\bibitem[Boyd and Mattingley, 2007]{boyd2007branch}
Boyd, S. and Mattingley, J. (2007).
\newblock Branch and bound methods.
\newblock {\em Notes for EE364b, Stanford University}, pages 2006--07.

\bibitem[Danna et~al., 2005]{danna2005exploring}
Danna, E., Rothberg, E., and Pape, C.~L. (2005).
\newblock Exploring relaxation induced neighborhoods to improve mip solutions.
\newblock {\em Mathematical Programming}, 102(1):71--90.

\bibitem[D’Ambrosio et~al., 2010]{d2010experiments}
D’Ambrosio, C., Frangioni, A., Liberti, L., and Lodi, A. (2010).
\newblock Experiments with a feasibility pump approach for nonconvex minlps.
\newblock In {\em International Symposium on Experimental Algorithms}, pages
  350--360. Springer.

\bibitem[Eckstein and Nediak, 2007]{eckstein2007pivot}
Eckstein, J. and Nediak, M. (2007).
\newblock Pivot, cut, and dive: a heuristic for 0-1 mixed integer programming.
\newblock {\em Journal of Heuristics}, 13(5):471--503.

\bibitem[Fischetti et~al., 2005]{fischetti2005feasibility}
Fischetti, M., Glover, F., and Lodi, A. (2005).
\newblock The feasibility pump.
\newblock {\em Mathematical Programming}, 104(1):91--104.

\bibitem[Fischetti and Lodi, 2003]{fischetti2003local}
Fischetti, M. and Lodi, A. (2003).
\newblock Local branching.
\newblock {\em Mathematical programming}, 98(1):23--47.

\bibitem[Fischetti and Lodi, 2008]{fischetti2008repairing}
Fischetti, M. and Lodi, A. (2008).
\newblock Repairing mip infeasibility through local branching.
\newblock {\em Computers \& operations research}, 35(5):1436--1445.

\bibitem[Fischetti and Salvagnin, 2009]{fischetti2009feasibility}
Fischetti, M. and Salvagnin, D. (2009).
\newblock Feasibility pump 2.0.
\newblock {\em Mathematical Programming Computation}, 1(2):201--222.

\bibitem[Ghosh, 2007]{ghosh2007dins}
Ghosh, S. (2007).
\newblock Dins, a mip improvement heuristic.
\newblock In {\em International Conference on Integer Programming and
  Combinatorial Optimization}, pages 310--323. Springer.

\bibitem[Gleixner et~al., 2021]{miplib}
Gleixner, A., Hendel, G., Gamrath, G., Achterberg, T., Bastubbe, M., Berthold,
  T., Christophel, P.~M., Jarck, K., Koch, T., Linderoth, J., L\"ubbecke, M.,
  Mittelmann, H.~D., Ozyurt, D., Ralphs, T.~K., Salvagnin, D., and Shinano, Y.
  (2021).
\newblock {MIPLIB 2017: Data-Driven Compilation of the 6th Mixed-Integer
  Programming Library}.
\newblock {\em Mathematical Programming Computation}.

\bibitem[Glover and Laguna, 1997]{glover1997tabu}
Glover, F. and Laguna, M. (1997).
\newblock Tabu search kluwer academic.
\newblock {\em Boston, Texas}.

\bibitem[{Gurobi Optimization, LLC}, 2022]{gurobi}
{Gurobi Optimization, LLC} (2022).
\newblock {Gurobi Optimizer Reference Manual}.

\bibitem[Hanafi and Todosijevi{\'c}, 2017]{hanafi2017mathematical}
Hanafi, S. and Todosijevi{\'c}, R. (2017).
\newblock Mathematical programming based heuristics for the 0--1 mip: a survey.
\newblock {\em Journal of Heuristics}, 23(4):165--206.

\bibitem[Hansen et~al., 2006]{hansen2006variable}
Hansen, P., Mladenovi{\'c}, N., and Uro{\v{s}}evi{\'c}, D. (2006).
\newblock Variable neighborhood search and local branching.
\newblock {\em Computers \& Operations Research}, 33(10):3034--3045.

\bibitem[Hung et~al., 1999]{hung1999evolutionary}
Hung, Y., Shih, C., and Chen, C. (1999).
\newblock Evolutionary algorithms for production planning problems with setup
  decisions.
\newblock {\em Journal of the Operational Research Society}, 50(8):857--866.

\bibitem[Lawler and Wood, 1966]{lawler1966branch}
Lawler, E.~L. and Wood, D.~E. (1966).
\newblock Branch-and-bound methods: A survey.
\newblock {\em Operations research}, 14(4):699--719.

\bibitem[Lokketangen and Glover, 1998]{lokketangen1998solving}
Lokketangen, A. and Glover, F. (1998).
\newblock Solving zero-one mixed integer programming problems using tabu
  search.
\newblock {\em European journal of operational research}, 106(2-3):624--658.

\bibitem[Medaglia et~al., 2009]{medaglia2009hybrid}
Medaglia, A.~L., Villegas, J.~G., and Rodr{\'\i}guez-Coca, D.~M. (2009).
\newblock Hybrid biobjective evolutionary algorithms for the design of a
  hospital waste management network.
\newblock {\em Journal of Heuristics}, 15(2):153--176.

\bibitem[Mitten, 1970]{mitten1970branch}
Mitten, L. (1970).
\newblock Branch-and-bound methods: General formulation and properties.
\newblock {\em Operations Research}, 18(1):24--34.

\bibitem[Momeni and Sarmadi, 2016]{momeni2016genetic}
Momeni, M. and Sarmadi, M. (2016).
\newblock A genetic algorithm based on relaxation induced neighborhood search
  in a local branching framework for capacitated multicommodity network design.
\newblock {\em Networks and Spatial Economics}, 16(2):447--468.

\bibitem[Nair et~al., 2020]{nair2020solving}
Nair, V., Bartunov, S., Gimeno, F., von Glehn, I., Lichocki, P., Lobov, I.,
  O'Donoghue, B., Sonnerat, N., Tjandraatmadja, C., Wang, P., et~al. (2020).
\newblock Solving mixed integer programs using neural networks.
\newblock {\em arXiv preprint arXiv:2012.13349}.

\bibitem[Rothberg, 2007]{rothberg2007evolutionary}
Rothberg, E. (2007).
\newblock An evolutionary algorithm for polishing mixed integer programming
  solutions.
\newblock {\em INFORMS Journal on Computing}, 19(4):534--541.

\bibitem[Schrijver et~al., 2003]{schrijver2003combinatorial}
Schrijver, A. et~al. (2003).
\newblock {\em Combinatorial optimization: polyhedra and efficiency},
  volume~24.
\newblock Springer.

\end{thebibliography}

\section*{Appendix}
We report in the following tables the computational results obtained in each of the 800 test instances of MIPLIB library we have used to evaluate our PUSH against Feasbility Pump. Since the success rate of the three rounding algorithms is significantly lower than the one of Feasibility and PUSH, those results are not reported. The symbol $-$ indicates  the algorithm was not able to find a feasible solution after 250 iterations.
\begin{table}[!ht]
\centering
\tiny

\begin{tabular}
    {l l l l l l l l l l} \hline
    
       Instance&\multicolumn{2}{c}{FP} &  \multicolumn{2}{c}{PUSH}& $N_{vars}$ & $N_{int}$ & $N_{bin}$  &$N_{constr}$ & $N_{nonzeros}$\\ \hline
    
         & Time& Gap & Time & Gap&  \\ \cline{2-5}
        ej & 291.54 & 396.99 & 900.02 &  -   & 3 & 0 & 3 & 1 & 3 \\ 
        flugpl & 96.43 & 0.21 & 900.01 &  -   & 18 & 0 & 11 & 18 & 46 \\ 
        stein9inf & 900 &  -   & 900 &  -   & 9 & 9 & 0 & 14 & 54 \\ 
        flugplinf & 900.04 &  -   & 900.01 &  -   & 18 & 0 & 11 & 19 & 64 \\ 
        enlight4 & 900.03 &  -   & 900.01 &  -   & 32 & 16 & 16 & 16 & 80 \\ 
        gr4x6 & 0.02 & 0.37 & 0.02 & 0.48 & 48 & 24 & 0 & 34 & 96 \\
\hline
        \end{tabular}
\end{table}

\begin{table}[!ht]
\centering
\tiny

    \begin{tabular}
    {l l l l l l l l l l} \hline
    
       Instance&\multicolumn{2}{c}{FP} &  \multicolumn{2}{c}{PUSH}& $N_{vars}$ & $N_{int}$ & $N_{bin}$  &$N_{constr}$ & $N_{nonzeros}$\\ \hline
    
         & Time& Gap &  Time & Gap&  \\ \cline{2-5}   milo-v13-4-3d-4-0 & 1637.54 &  -   & 900.05 &  -   & 688 & 160 & 0 & 1328 & 3292 \\ 
        22433 & 2.44 & 0 & 0.12 & 0.01 & 429 & 231 & 0 & 198 & 3408 \\ 
        neos-911970 & 0.06 & 2.76 & 0.09 & 2.76 & 888 & 840 & 0 & 107 & 3408 \\ 
        qiu & 0.5 & 8.03 & 0.34 & 7.17 & 840 & 48 & 0 & 1192 & 3432 \\ 
        haprp & 1.36 & 0.01 & 1.91 & 0 & 1828 & 0 & 1828 & 1048 & 3628 \\ 
        neos-3421095-cinca & 0.07 & 0.17 & 0.06 & 0.17 & 896 & 117 & 42 & 1218 & 3678 \\
\hline
  
                \end{tabular}
\end{table}

\begin{table}[!ht]
\centering
\tiny

    \begin{tabular}
    {l l l l l l l l l l} \hline
    
       Instance&\multicolumn{2}{c}{FP} &  \multicolumn{2}{c}{PUSH}& $N_{vars}$ & $N_{int}$ & $N_{bin}$  &$N_{constr}$ & $N_{nonzeros}$\\\hline 
        & Time& Gap &  Time & Gap&  \\ \cline{2-5}        p500x2988 & 900.1 &  -   & 719.95 & 308.27 & 5976 & 2988 & 0 & 3488 & 11952 \\ 
        p500x2988c & 1.63 & 60.96 & 1.23 & 24.34 & 5976 & 2988 & 0 & 3488 & 11952 \\ 
        p500x2988d & 1.14 & 108.98 & 0.73 & 59.66 & 5976 & 2988 & 0 & 3488 & 11952 \\ 
        10teams & 10.19 & 0.05 & 1.93 & 0.06 & 2025 & 1800 & 0 & 230 & 12150 \\ 
        h50x2450 & 1.29 & 3.81 & 4.82 & 3.72 & 4900 & 2450 & 0 & 2549 & 12152 \\ 
        supportcase30 & 900.29 &  -   & 900.19 &  -   & 1024 & 1024 & 0 & 1028 & 12288 \\
\hline
 
    \end{tabular}
\end{table}

\begin{table}[!ht]
\centering
\tiny

    \begin{tabular}
    {l l l l l l l l l l} \hline
    
       Instance&\multicolumn{2}{c}{FP} &  \multicolumn{2}{c}{PUSH}& $N_{vars}$ & $N_{int}$ & $N_{bin}$  &$N_{constr}$ & $N_{nonzeros}$\\ \hline
       & Time& Gap &  Time & Gap& \\ \cline{2-5}              neos-4393408-tinui & 900.13 &  -   & 900.12 &  -   & 4004 & 2000 & 0 & 4554 & 23496 \\ 
        sct2 & 0.56 & 1.35 & 0.52 & 1.49 & 5885 & 2872 & 0 & 2151 & 23643 \\ 
        bppc4-08 & 0.17 & 1.02 & 0.35 & 1.25 & 1456 & 1454 & 0 & 111 & 23964 \\ 
        neos-3083819-nubu & 51.88 & 0.04 & 20.66 & 0.03 & 8644 & 0 & 8644 & 4725 & 24048 \\ 
        neos18 & 0.45 & 0.19 & 3.33 & 0.31 & 3312 & 3312 & 0 & 11402 & 24614 \\ 
        neos-1582420 & 900.16 &  -   & 900.16 &  -   & 10100 & 10000 & 100 & 10180 & 24814 \\
\hline
  
    \end{tabular}
\end{table}

\begin{table}[!ht]
\centering
\tiny

    \begin{tabular}
    {l l l l l l l l l l} \hline
    
       Instance&\multicolumn{2}{c}{FP} &  \multicolumn{2}{c}{PUSH}& $N_{vars}$ & $N_{int}$ & $N_{bin}$  &$N_{constr}$ & $N_{nonzeros}$\\ \hline
       & Time& Gap &  Time & Gap& \\ \cline{2-5}

        pythago7825 & 900.76 &  -   & 900.08 &  -   & 3745 & 3745 & 0 & 14672 & 44016 \\ 
        eil33-2 & 8.85 & 1.53 & 145.05 & 0.97 & 4516 & 4516 & 0 & 32 & 44243 \\ 
        neos-480878 & 0.96 & 0.59 & 1.36 & 0.47 & 534 & 189 & 0 & 1321 & 44370 \\ 
        neos-3045796-mogo & 900.15 &  -   & 900.16 &  -   & 11016 & 11016 & 0 & 2226 & 44442 \\ 
        s1234 & 900.08 &  -   & 901.17 &  -   & 2945 & 2945 & 0 & 8418 & 44641 \\ 
        CMS750\_4 & 900.18 &  -   & 704.73 & 2.21 & 11697 & 7196 & 0 & 16381 & 44903 \\
\hline
 
    \end{tabular}
\end{table}

\begin{table}[!ht]
\centering
\tiny

    \begin{tabular}
    {l l l l l l l l l l} \hline
    
       Instance&\multicolumn{2}{c}{FP} &  \multicolumn{2}{c}{PUSH}& $N_{vars}$ & $N_{int}$ & $N_{bin}$  &$N_{constr}$ & $N_{nonzeros}$\\ \hline
        & Time& Gap &  Time & Gap& \\ \cline{2-5}
        
        wachplan & 14.83 & 0 & 3.27 & 0 & 3361 & 3360 & 1 & 1553 & 89361 \\ 
        neos-952987 & 900.3 &  -   & 900.35 &  -   & 31329 & 31329 & 0 & 354 & 90384 \\ 
        air03 & 1.42 & 0.09 & 5.03 & 0.17 & 10757 & 10757 & 0 & 124 & 91028 \\ 
        neos-3216931-puriri & 900.13 &  -   & 900.15 &  -   & 3555 & 3268 & 0 & 5989 & 91691 \\ 
        msc98-ip & 900.25 &  -   & 900.27 &  -   & 21143 & 20237 & 53 & 15850 & 92918 \\ 
        neos-3581454-haast & 900.2 &  -   & 900.21 &  -   & 8112 & 7512 & 0 & 17220 & 92976 \\
\hline

    \end{tabular}
\end{table}

\begin{table}[!ht]
\centering
\tiny

    \begin{tabular}
    {l l l l l l l l l l} \hline
    
       Instance&\multicolumn{2}{c}{FP} &  \multicolumn{2}{c}{PUSH}& $N_{vars}$ & $N_{int}$ & $N_{bin}$  &$N_{constr}$ & $N_{nonzeros}$\\ \hline
        & Time& Gap &  Time & Gap& \\ \cline{2-5}  
        neos-953928 & 2.06 & 0.92 & 8.25 & 0.31 & 23305 & 23305 & 0 & 12498 & 169861 \\ 
        stockholm & 902.1 &  -   & 902.11 &  -   & 20644 & 962 & 0 & 57346 & 171076 \\ 
        neos-5266653-tugela & 900.3 &  -   & 900.29 &  -   & 23310 & 22662 & 0 & 23458 & 171534 \\ 
        neos-5196530-nuhaka & 900.24 &  -   & 900.26 &  -   & 14544 & 9792 & 0 & 42124 & 171702 \\ 
        ger50-17-trans-dfn-3t & 900.75 &  -   & 612.63 &  -   & 22414 & 0 & 18062 & 499 & 172035 \\ 
        ger50-17-trans-pop-3t & 900.36 &  -   & 900.3 &  -   & 22414 & 0 & 18062 & 499 & 172035 \\
\hline

        \end{tabular}
\end{table}

\begin{table}[!ht]
\centering
\tiny

    \begin{tabular}
    {l l l l l l l l l l} \hline
    
       Instance&\multicolumn{2}{c}{FP} &  \multicolumn{2}{c}{PUSH}& $N_{vars}$ & $N_{int}$ & $N_{bin}$  &$N_{constr}$ & $N_{nonzeros}$\\ \hline 
        & Time& Gap &  Time & Gap&  \\ \cline{2-5}  
        gasprod2-1 & 900.41 &  -   & 900.73 &  -   & 29568 & 7389 & 0 & 73448 & 282846 \\ 
        satellites2-25 & 900.3 &  -   & 900.19 &  -   & 35378 & 34324 & 0 & 20916 & 283668 \\ 
        satellites2-40 & 905.36 &  -   & 900.19 &  -   & 35378 & 34324 & 0 & 20916 & 283668 \\ 
        sp97ar & 4.09 & 29 & 5.07 & 7.56 & 14101 & 14101 & 0 & 1761 & 290968 \\ 
        bts4-cta & 900.5 &  -   & 901.31 &  -   & 75400 & 2260 & 0 & 45350 & 291904 \\ 
        cryptanalysiskb128n5obj14 & 900.48 &  -   & 900.32 &  -   & 48950 & 47830 & 1120 & 98021 & 292875 \\
\hline
        \end{tabular} 
\end{table}

\end{document}